\documentclass [12pt] {article}
\usepackage{amsfonts}
\usepackage{amssymb}
\newcommand{\qed} {\hspace {0.1in} \rule {1.5mm} {3.5mm}}

\newtheorem{lemma}{Lemma}[section]
\newtheorem{corollary}[lemma]{Corollary}
\newtheorem{theorem}{Theorem}

\newtheorem{conjecture}{Conjecture}
\newtheorem{proposition}[lemma]{Proposition}
\newtheorem{definition}[lemma]{Definition}
\def\fix  {{\#_{\rm fix}}}
\def\lift#1{{\widetilde{#1}}}
\def\Sn   {{\mathfrak S(\tau)}}
\def\Sno  {{\mathfrak S(\tau/\omega)}}

\def\proof {\smallskip\noindent{\bf Proof:} }

\newcommand{\Mat}[0]{{\rm Mat}}
\newcommand{\Proj}[0]{{\rm Proj}}
\newcommand{\cN}[0]{{\mathcal N}}
\newcommand{\Fin}[0]{{\mathcal Fin}}
\newcommand{\Gclass}[0]{{\mathcal G}}
\newcommand{\Fix}[0]{{\rm Fix}}
\newcommand{\Id}[0]{{\rm Id}}
\newcommand{\Tr}[0]{{\rm Tr}}
\newcommand{\tra}[0]{{\rm tr}}
\newcommand{\C}[0]{{\mathbb C}}
\newcommand{\bN}[0]{{\mathbb N}}
\newcommand{\bR}[0]{{\mathbb R}}
\newcommand{\bZ}[0]{{\mathbb Z}}
\newcommand{\Inv}[0]{{\rm Inv}}

\title{Hyperlinearity, essentially free actions and $L^2$-invariants. The
sofic property}

\author{{\sc G\'abor Elek and Endre Szab\'o}
\cr Mathematical Institute of
the Hungarian Academy of Sciences\cr P.O. Box 127, H-1364 Budapest, Hungary\cr
elek@renyi.hu, endre@renyi.hu}
\date{}
\begin{document}

\maketitle
\noindent{\bf Abstract.} We prove that Connes' Embedding Conjecture holds for
the von Neumann algebras of sofic groups, that is sofic groups are
hyperlinear. Hence we provide some new examples of hyperlinearity. We also
show that the Determinant Conjecture holds for sofic groups as well. We
introduce the notion of essentially free actions and amenable actions and
study their properties.
\vskip 0.2in
\newpage

%%%    1. oldal

\section{Introduction and preliminaries}

\subsection{The Gromov-Hausdorff topology}

A finitely generated marked group $(\Gamma,S)$ is a group $\Gamma$ with an ordered
generator system $S=\{s_1, s_2,\dots, s_m\}$. The set of such marked groups
of $m$ generators is denoted by $\Sigma_m$. The Gromov-Hausdorff metric on
$\Sigma_m$
is given as follows \cite{Champ}:
$$d_{\Sigma_m}((\Gamma,S),(\Gamma',S'))=2^{-R}$$ if the two groups have exactly the same
relations of length at most $R$ with respect to their marked generators, on
the other hand a relation of length $R+1$ holds only for one of them. Then
$\Sigma_m$ is a compact metric space.
Let $B_{(\Gamma,S)}(k)$ be the ball of radius $k$ in the Cayley-graph of $(\Gamma,S)$
labeled by$\{1,2,\dots,m\}$, where $(x\to y)$ is a directed edge of the
Cayley-graph labeled by $i$ if $s_i x=y$. If $d_{\Sigma_m}((\Gamma,S),(\Gamma',S'))\leq
2^{-(2k+1)}$ then $B_{(\Gamma,S)}(k)$ and $B_{(\Gamma',S')}(k)$ are isometric as labeled
graphs.

Obviously, if $\{(\Gamma_n,S_n)\}_{n\geq 1}$ is a convergent sequence of
torsion-free marked groups, then $\lim_{n\to\infty}\{(\Gamma_n,S_n)\}
=(\Gamma,S)$ is torsion-free as
well.

\subsection{Sofic groups}
In \cite{VG} the authors studied the limit groups of finite marked groups,
which they called $LEF$-groups. This lead Gromov to the notion of sofic groups
\cite{Gro},\cite{Wei}.
\begin{definition} \label{sofic1}
Let $\Gamma$ be a finitely generated group and $S\subseteq \Gamma$ be a finite symmetric
generator set. Then the group $\Gamma$ is called \it{sofic} if there exists a
sequence of finite directed graphs $\{(V_n,E_n)\}_{n\geq 1}$, edge-labeled by
$S$ and subsets $V_n^0\subseteq V_n$ with the following property:

\noindent
For any $\delta>0$ and $r\in \bN$, there is an integer $n_{r,\delta}>0$
such that if $m\geq n_{r,\delta}>0$, then

\begin{itemize}
\item For each $v\in V^0_m$, there is a map $\phi_v: B_{(\Gamma,S)}(r)\to
V_m$, which is an isomorphism (of labeled graphs) between $B_{(\Gamma,S)}(r)$
and the $r$-ball in $V_m$ around $v$.

\item $|V^0_m|\geq (1-\delta)|V_m|\,.$ \end{itemize}
\end{definition}
Note that we may suppose that the graphs $\{(V_n,E_n)\}_{n\geq 1}$ have uniformly
bounded vertex degrees.

\noindent
We call a not necessarily finitely generated group $\Gamma$ sofic 
if all of its finitely generated subgroups are sofic. However one has an
equivalent definition for soficity as well \cite{ESZ2}:
\begin{definition} \label{sofic-group}
The group $\Gamma$ is sofic if for any real number $0<\epsilon<1$ and any finite
subset $F\subseteq \Gamma$ there exists a natural number $n$ and a function
$\psi_n:\Gamma\to S_n$ from $\Gamma$ into the group of permutations on $n$ elements with
the following properties:
\begin{enumerate}
\item[(a)]   $\fix\Big(\phi(e)\phi(f)\,\phi(ef)^{-1}\Big)\ge(1-\varepsilon)n$
             \kern20pt for any two elements $e,f\in F$.
\item[(b)]   $\phi(1)=1$.
\item[(c)]   $\fix\phi(e)\le\varepsilon n$\kern20pt
             for any $1\ne e\in F$,
\end{enumerate}
where $\fix\pi$ denotes the number of fixed points of the permutation $\pi\in S_n$.
\end{definition}
It is known \cite{Gro},\cite{Wei},\cite {ESZ1} that the $LEF$-groups are
sofic. Direct products, subgroups, free products, inverse and direct limits of
sofic groups are sofic as well. If $N\lhd \Gamma$, $N$ is sofic and $\Gamma/N$ is amenable, then $\Gamma$ is also
sofic.
Residually amenable groups are sofic, however there exist finitely generated 
non-residually amenable sofic groups as well \cite{ESZ1}. It is conjectured that there are
non-sofic groups, but no example is known yet. In this paper we prove that any
countable sofic group can be imbedded into a  countable simple sofic group,
particularly, there exist countable simple non-amenable sofic groups (Theorem \ref{generic-sofic}).

\subsection{Connes' Embedding Conjecture and hyperlinear groups}
Let $\mathcal H$ be a separable Hilbert-space with orthonormal basis
$\{e_n\}_{n=1}^\infty$. We shall denote by $B(\mathcal H)$ the space
of bounded linear operators on $\mathcal H$. Then for any $k\ge1$, we
have an embedding of the matrix ring
$\rho_k:\Mat_{k\times k}(\C)\to B(\mathcal H)$
given by the formula
$$
\rho_k(A)(f_j^s)=
\sum_{i=1}^{k}A_{ij}f_i^s
\kern 25pt\mbox{with \ }f_i^s=e_{ks+i}
$$
for each matrix $A=\{a_{ij}\}_{1\leq i,j\leq k}\in\Mat_{k\times k}(\C)$,
all indices $1\le j\le k$ and $s\ge 0$.
It is easy to see that if $l$ is a multiple of $k$ then the image of
$\rho_l$ contains the image of $\rho_k$. In fact there is a natural
inclusion
$T_{kl}:\Mat_{k\times k}\to
\Mat_{k\times k}(\C)\otimes\Mat_{l/k\times l/k}(\C)\cong
\Mat_{l\times l}(\C)$
given by the formula $a\to a\otimes\Id$
(here $\Id$ is the $l/k \times l/k$ unit matrix), and all of these embeddings
are compatible:
$\rho_l\circ\mathcal T_{kl}=\rho_k$.
Hence the images of the $\rho_k$ form a directed system of subrings in
$B(\mathcal H)$.
The {\em hyperfinite factor} $R$
 is the closure of this union
in the ultraweak operator topology.
%%
%%   2. oldal
%%
The natural normalized traces on $\Mat_{k\times k}(\C)$, $k\geq 1$ define a
normalized trace on the union of their images, and this extends onto
$R$ as a tracial state $\Tr_R$.

Now let $\omega$ be an ultrafilter on a countable set $I$,
let $\lim_\omega$ denote the corresponding ultralimit.
Consider the restricted direct product:
$$
R^{(I)} =
\Bigg\{\{a_i\}_{i\in I}\in\prod_{i\in I} R\ \Bigg|\
\sup_{i\in I}\|a_i\|<\infty\Bigg\}\,.
$$
Then
$$
J_\omega =
\Bigg\{\{a_i\}_{i\in I}\in\prod_{i\in I} R\ \Bigg|\
\lim_\omega\Tr_R(a_ia_i^*)=0\Bigg\}
$$
forms a maximal ideal and
$$
R^\omega=R^{(I)}/J_\omega
$$
is a type II$_1$-factor with trace
$$
\Tr_\omega\Big(\{a_i\}_{i\in I}\Big) =
\lim_\omega\Tr_R(a_i).
$$
Connes conjectured that any type II$_1$-factor of separable predual
can be embedded into $R^\omega$.

%% 3. oldal

It is known that if a countable group can be faithfully represented as
unitary elements of $R^\omega$ then it's von Neumann algebra can be
embedded into $R^\omega$. Such groups are called {\em hyperlinear}
\cite{Rad}, \cite{Oza}.
Note that groups with the so-called factorization property (see
\cite{Kirch} for a definition) are always hyperlinear. In \cite{Oza}
Ozawa showed that the  groups with factorization property
are closed when taking subgroups, extensions by amenable quotients,
free products, and increasing union (like sofic groups). He also proved that the
factorization property is residual. 

%%  4. oldal

We shall prove that countable sofic groups are  hyperlinear ( Theorem \ref{x}).
 Thus we obtain examples of simple non-amenable and finitely generated
 non-residually amenable
 hyperlinear groups.

\subsection{Amenable actions}

The notion of amenable action was introduced by Tarski. Suppose that a
group $\Gamma$ acts on a set $X$. Then the following conditions are
equivalent \cite{Pat}:

%% 5. oldal

\begin{enumerate}
\item	There exists a finitely additive $\Gamma$-invariant
probability measure on $X$.
\item	There is no paradoxical decomposition of $X$ with respect to
the $\Gamma$-action. That is, $X$ cannot be partitioned into subsets
$A_1,A_2,\dots A_l$ such that for some elements $g_1,h_1,g_2,h_2,\dots,
g_l,h_l\in\Gamma$ the subsets
$g_1A_1, h_1A_1, g_2A_2, h_2A_2,\dots, g_lA_l, h_lA_l$
are pairwise disjoint.
\item	There exist F\o lner-systems. That is, for any finite subset
$K\subseteq\Gamma$ and any $\epsilon>0$ one can find a finite subset
$F\subseteq X$ such that for any $g\in K$
$$
\frac{|gF\vartriangle F|}{|F|}<\epsilon.
$$
\end{enumerate}

%%  6. oldal

If an action satisfies the above conditions then it is called an
{\em amenable action}. If $\Gamma$ is an amenable group, then all its
actions are amenable, however any group acts on a single point in an
amenable way.

\begin{definition}   \label{essentially-freee-action-def}
An amenable action of a group $\Gamma$ on a set $X$ is called
{\em essentially free} if there exists a $\Gamma$-invariant finitely
additive measure $\mu$ on $X$ such that for any $1\ne g\in\Gamma$,
$$
\mu\Big(\Fix(g)\Big)=0,
$$
where $\Fix(g)\subseteq X$ denotes the fixed point set of $g$.
\end{definition}

We shall see that amenable groups and residually free groups
have essentially free amenable actions.
On the other hand if a group $\Gamma$ of Kazhdan's Property (T) has an
essentially free amenable action, then $\Gamma$ must be residually
finite. Although, we are unable to construct a non-residually finite
sofic group with Kazhdan's Property (T) we show that the non-existence
of such groups would imply the existence of a residually finite
word-hyperbolic group.

%%  7. oldal

We introduce the notion of an almost-action.

\begin{definition} \label{almost-action-def}
Let $\Gamma$ be a group, $X$ be a set and $\mu$ be a finitely additive
probability measure on $X$. For any $g\in\Gamma$, let $\phi(g):X\to X$
be a bijection. We say that $(\Gamma,X,\mu,\phi)$ is an {\em amenable
almost-action} if
\begin{enumerate}
\item	$\phi(1)$ is the identity transformation on $X$.
\item	For any $g\in\Gamma$, $\phi(g)$ is measure-preserving:
$\mu\big(\phi(g)(L)\big)=\mu(L)$ for any $L\subseteq X$.
\item	For any $g,h\in\Gamma$ there exists a subset
$A_{g,h}\subseteq X$ of measure $\mu(A_{g,h})=1$ such that
$\phi(gh)(x)=\phi(g)\big(\phi(h)(x)\big)$ if $x\in A_{g,h}$.
\end{enumerate}
\end{definition}

%%  8. oldal

\begin{definition}  \label{essentially-freee-almost-action-def}
Again, we call an amenable almost-action $(\Gamma,X,\mu,\phi)$
{\em essentially free}, if $\mu(\Fix(\phi(g)))=0$ for all
$g\in\Gamma$.
\end{definition}

We shall obtain Tarski-type characterizations for essentially free
amenable actions and almost-actions, and prove that a group
has an essentially free amenable almost-action if and only if
it is sofic.
Finally we prove that if $\Gamma$ is an inner amenable group then
either $\Gamma$ has a nontrivial sofic quotient or $\Gamma$ is almost
commutative (Corollary \ref{inner}).

\subsection{$L^2$-invariants}
Let $\Gamma$ be a countable group and $N\Gamma$ be the von Neumann algebra with its
usual trace $\Tr_\Gamma$. Then one can extend the trace to $\Mat_{d\times d}(N\Gamma)$
by
$$\Tr_\Gamma(A):=\sum^d_{i=1} \Tr_\Gamma(a_{ii})\,.$$
Let $\Delta$ be a positive self-adjoint operator in $\Mat_{d\times d}(N\Gamma)$,
then one can consider the spectral density function
$$F_\Delta(\lambda):=Tr_\Gamma\,\chi[0,\lambda] (\Delta)\,$$
$$F_\Delta(0):=Tr_\Gamma\,Proj(\Delta)=\dim_\Gamma(\ker \Delta)\,$$
where $\Proj(\Delta)$ is the orthogonal projection onto the kernel of
$\Delta$,
and $\dim_\Gamma$ is the von Neumann-dimension.

\noindent
The Fuglede-Kadison determinant of $\Delta$ is defined as follows.
\[ \ln\det_{\cN}(\Delta):=\left\{ \begin{array}{ll}
\int^\infty_{0^+}\,\ln(\lambda)\,dF_{\Delta}(\lambda) & 
\mbox{if the integral converges}\\
-\infty & \mbox{otherwise}
\end{array}
\right. \]
There are two major conjectures concerning the invariants above.

\noindent
\begin{conjecture}{\bf The Determinant Conjecture:}
If $\Delta\in\Mat_{d\times d}(\bZ \Gamma)$ is a positive self-adjoint operator,
then
$\ln\det_{\cN}(\Delta)\geq 0\,.$
\end{conjecture}
\begin{conjecture}{\bf The Atiyah Conjecture:}
If $\Delta\in\Mat_{d\times d}(\bZ \Gamma)$ is a positive self-adjoint operator and
$\Gamma$ is torsion-free, then $F_{\Delta}(0)$ is an integer.
\end{conjecture}
In \cite{Schick1}\cite{Schick2}, Schick proved the Determinant Conjecture
 for the large class
of groups $\Gclass$ and the Atiyah Conjecture for residually torsion-free elementary
amenable groups. In \cite{Schick1} the author remarked that no example of 
non-residually amenable group is known in the class $\Gclass$. We prove that
the Determinant Conjecture actually holds for sofic groups
 (Theorem \ref{soficdet}). Particularly, we
show the existence of finitely generated non-residually amenable groups and
simple non-amenable groups satisfying the conjecture. We also show that a more
geometric version of Schick's second result holds. That is if in $\Sigma_m$
the Atiyah Conjecture holds for a convergent sequence of torsion-free sofic
groups then it must hold for the limit group as well (Theorem \ref{soficlimit}).

%%%%%%%%%%%%%%%%%%%%%%%%%%%%%%%%%%%%%%%%55

\section{Universal sofic groups}\label{secsofic}
In this section we give a new characterization of sofic groups.
There is a family $\Sno$ of sofic groups with very explicit
definition, and each sofic group is isomorphic to a subgroup of one of
these $\Sno$.
\begin{definition}   \label{symmetric-group}
For an integer $d$ we denote by $S_d$ the full symmetric group on
the set $\{1,2,\dots d\}$ (we use left actions).
For a permutation $\alpha\in S_d$, let $\#_t\alpha$ denote the
number of $t$-cycles in $\alpha$, and $\fix\alpha$ denotes the
number of fixed points of $\alpha$.
Clearly $\fix\alpha=\#_1\alpha$.
For a natural number $k\in\mathbb N$ we can identify
the $k$-fold disjoint union $\coprod^k\{1,2,\dots, d\}$
with $\{1,2,\dots, kd\}$ and the $k$-fold direct product
$\prod^k\{1,2,\dots, d\}$ with $\{1,2,\dots, d^k\}$.
This induces natural homomorphisms
$\coprod^k:S_d\to S_{kd}$ and $\prod^k:S_d\to S_{d^k}$.
It is clear that
$$
\fix(\coprod^k\alpha)  =  k\cdot\fix\alpha
\kern 20pt\mbox{and}\kern20pt
\fix(\prod^k\alpha)    =  (\fix\alpha)^k
$$
for all $k\in\mathbb N$ and  $\alpha\in S_d$.
\end{definition}

\begin{definition}
\label{sno-def}
Let $\omega$ be an ultrafilter on a set $I$, and $\tau:I\to\mathbb N$ a
function with the property $\lim_\omega \tau=\infty$.
Let $\Sn$ denote the direct product  $\prod_{i\in I}S_{\tau(i)}$,
and for elements $\pi\in\Sn$ and indices $i\in I$ let $\pi(i)\in S_{\tau(i)}$
denote the $i$-th component of $\pi$. We define the average number of
$t$-cycles and the average number of fixed points of $\pi$ via the formula:
$$
\#_t\pi = \lim_\omega\frac{\#_t\pi(i)}{\tau(i)}
\kern 40pt
\fix\pi = \lim_\omega\frac{\fix\pi(i)}{\tau(i)}\,.
$$
Then we define the quotient group:
$$
\Sno = \Sn \Bigg/
	\left\{ \pi\in \Sn \Big|\fix\pi = 1\right\}\,.
$$
For elements $\pi\in\Sno$ we choose a coset representative
$\lift\pi\in\Sn$, and define the average number of $t$-cycles and
fixpoints as $\#_t\pi=\#_t\lift\pi$ and
$\fix\pi=\fix\lift\pi$.
\end{definition}

\begin{proposition} \label{sno-properties}
Let $I$, $\omega$ and $\tau$ as in Definition~\ref{sno-def}.
\begin{enumerate}
\item	\label{well-defined-sno}
	The  subgroup  of $\Sn$ in
	Definition~\ref{sno-def} is indeed a normal
	subgroup, so $\Sno$ is well-defined.
\item	\label{well-defined-cycle-number}
	The average number of $t$-cycles and fixed points of an element
	$\pi\in\Sno$ does not depend  on the chosen representative
	$\lift\pi$. So $\#_t\pi$ and $\fix\pi$ is well-defined on
	$\Sno$.
\item	\label{cycle-number-properties}
	The function $\#_t$ is a conjugation invariant function on
	$\Sno$, and $\sum_{t=1}^\infty t\cdot\#_t\pi\le 1$ for each
	element $\pi\in\Sno$. Moreover, if $P_t$ is a sequence of
	nonnegative real numbers satisfying
	$\sum_{t=1}^\infty t\cdot P_t\le 1$, then there is an element
	$\pi\in\Sno$ which satisfies  $\#_t\pi = P_t$ for all $t$.
\item	\label{conjugacy-classes-in-sno}
	The conjugacy classes of $\Sno$ are uniquely determined by
	their $\#_t$ invariants: elements $\pi,\rho\in\Sno$ are
conjugate if and only if $\#_t\pi=\#_t\rho$ for all $t$.
\item	\label{simple-sno}
	Each nontrivial conjugacy class generates the group $\Sno$.
	Hence it is a simple group.
\end{enumerate}
\end{proposition}

{\bf Proof of}
(\ref{well-defined-sno}) and (\ref{well-defined-cycle-number}):
The number of $t$-cycles is a conjugation invariant function on the
symmetric group $S_d$,
hence $\#_t$ and $\fix$ are conjugation invariant functions on
$\Sn$. For permutations $\alpha,\beta\in S_d$ those $t$-cycles of
$\alpha$ which lie entirely in the fixed point set of $\beta$ will also
be $t$-cycles of the composition $\alpha\beta$, so
$\#_t(\alpha\cdot\beta)\ge\#_t\alpha+\fix\beta-d$.
Hence $\#_t(\pi\rho)\ge\#_t\pi+\fix\rho-1$
for all $\pi,\rho\in\Sn$.
This implies that the subset $N=\{\rho\in\Sn | \fix\rho=1\}$ is indeed a
normal subgroup of $\Sn$, and for all $\rho\in N$ and all
$\pi\in\Sn$ one has $\#_t(\pi\rho)=\#_t\pi$. Hence the
average number of $t$-cycles and the average number of fixed points is
in fact well defined on $\Sn/N=\Sno$.
\qed

{\bf Proof of}
(\ref{cycle-number-properties}):
Since $\#_t$ and $\fix$ are conjugacy invariant functions on the
symmetric groups $S_{\tau(i)}$, they are also conjugation invariant on
$\Sn$ and $\Sno$.
It is clear that for $\alpha\in S_d$ one has
$\sum_{t=1}^\infty t\cdot\#_t\alpha=d$, and this implies that
$\sum_{t=1}^\infty t\cdot\#_t\pi\le1$ for all $\pi\in\Sno$.
(Note that $\lim_\omega$ is finitely additive and monotonic, but not
infinitely additive.) Now let $P_t$ be  a sequence of nonnegative real
numbers such that $\sum_{t=1}^\infty t\cdot P_t\le1$. For each
index $i\in I$ we build a permutation $\pi_i\in S_{\tau(i)}$ which
consists of $[P_t\cdot \tau(i)]$ disjoint $t$-cycles for each $t$
(we denote by $[x]$ the integer part of the real number $x$),
and one extra cycle on the remaining elements (if there are any).
We set $\lift\pi=\prod_{i\in I}\pi_i\in\Sn$. Let $\pi$ denote
the image of $\lift\pi$ in $\Sno$. Since the number of $t$-cycles of
each $\pi_i$ is between $P_t\cdot \tau(i)-1$ and $P_t\cdot \tau(i)+1$, and
$\lim_\omega 1/\tau(i) =0$, we find that
$\#_t\pi= P_t$ for all $t$.
\qed

{\bf Proof of}
(\ref{conjugacy-classes-in-sno}):
Let $\pi,\rho\in\Sno$ be elements with $\#_t\pi=\#_t\rho$
for all $t$, we shall prove that they are conjugate to each other.
First we pick representatives $\lift\pi,\lift\rho\in\Sn$. For each
index $i\in I$ and each natural number $t$ let
$p(t,i)=\min(\#_t\lift\pi(i),\#_t\lift\rho(i))$, then
\begin{equation} \label{number-of-t-cycles}
\#_t\pi=\#_t\rho=\lim_\omega p(t,i)
\kern 40pt\mbox{for all }t\,.
\end{equation}
For each index $i\in I$ we divide the cycles of $\lift\pi(i)$ into two
groups. Let $P(i)$ be an arbitrary collection of cycles which contains exactly
$p(t,i)$ among the $t$-cycles for all $t$,
and let $E(i)$ denote the collection of those cycles which are not
included in $P(i)$. We shall prove that
\begin{equation}  \label{number-of-exceptional-points}
\lim_\omega\frac{|E(i)|}{\tau(i)} = 0\,.
\end{equation}
Indeed, let $e(t,i)$ denote the number of $t$-cycles in $E(i)$. Then
%$
%\lim_\omega e(t,i)/\tau(i) =
%\lim_\omega(\#_t\lift\pi(i)-p(t,i))/\tau(i)=0
%$
$$
\lim_\omega\frac{e(t,i)}{\tau(i)} =
\lim_\omega\frac{\#_t\lift\pi(i)-p(t,i)}{\tau(i)}=0
$$
by equation~(\ref{number-of-t-cycles}), and
$\sum_{t\ge T}e(t,i)\le \tau(i)/T$. Putting these together we see that
$$
\lim_\omega\sum_{t=1}^\infty\frac{e(t,i)}{\tau(i)} =
\sum_{t<T}\lim_\omega\frac{e(t,i)}{\tau(i)} +
\lim_\omega\sum_{t\ge T}\frac{e(t,i)}{\tau(i)}  \le
\frac{1}{T}
$$
for all $T>0$. This
proves equation~(\ref{number-of-exceptional-points}).
We pick one element from the support of each cycle in $E(i)$, let
$\varepsilon_i\in S_{\tau(i)}$ be an $|E(i)|$-cycle on these elements, and
let $\lift\varepsilon\in\Sn$ be the direct product of these $\varepsilon_i$.
Equation~(\ref{number-of-exceptional-points}) shows us that
$\fix\lift\varepsilon=1$, hence the element
$\lift\pi'=\lift\varepsilon\lift\pi\in\Sn$ is another representative of
$\pi$. It is clear from the construction that the cycles of
$\lift\pi'(i)$ are the cycles in $P(i)$, and possibly one extra cycle
on the remaining $\tau(i)-\sum_tp(t,i)$ elements of $\{1,2,\dots \tau(i)\}$.
(In fact, if there are no elements left, then this extra cycle is
missing.)
Similarly, we can replace $\lift\rho$ with another representative
$\lift\rho'\in\Sno$ with the property that it has $p(t,i)$ $t$-cycles
for each $t$, and possibly one extra cycle on the remaining
$\tau(i)-\sum_tp(t,i)$ elements.
But then $\lift\pi'(i)$ and $\lift\rho'(i)$ have the same number of
$t$-cycles for each $t$, so they are conjugate in $S_{\tau(i)}$. Hence
$\lift\pi'$ and $\lift\rho'$ are conjugate in $\Sn$, therefore
$\pi$ and $\rho$ are conjugate in $\Sno$.
This proves (\ref{conjugacy-classes-in-sno}).
\qed

\noindent
{\bf Proof of}
(\ref{simple-sno}):
\begin{lemma} \label{lvege1}
Let $\sigma$ be an element of $A_n$, the alternating group on $n$
elements. Suppose that $\sigma$ has an orbit of length two and
$n-2r\geq -1$, where $r$ is
the number of orbits of $\sigma$. Then $C^4_\sigma=A_n$, where $C_\sigma$
is the conjugacy class of $\sigma$.
\end{lemma}
{\bf Proof:}
This is just Theorem 3.05 in \cite{Brenner}. \qed
\begin{lemma} \label{lvege2}
Let $\sigma\in A_n$. Suppose that $\sigma$ has an orbit of length two. Let $T=n-\fix(\sigma)$. Then
$C^{8[\frac{n}{T}]}_\sigma=A_n.$
\end{lemma}
{\bf Proof:}
Suppose that $\frac{n}{2}\leq T < n$. Then by Lemma \ref{lvege1},
$C^4_\sigma$ contains an element $\sigma_1$ which moves exactly
$n-T$ elements and fixes the complete fixed point set of $\sigma$.
Therefore, for some $a\in A_n$, $\sigma a \sigma_1 a^{-1}$ has no
fixed points at all. Hence $C^8_\sigma=A_n$.

\noindent
If $T<\frac{n}{2}$ then there exist $[\frac{n}{T}]$ conjugates
of $\sigma$ such that they move mutually disjoint subsets.
Thus their product has less fixed points then $\frac{n}{2}$.
Consequently, by the previous argument, $C^{8[\frac{n}{T}]}_\sigma=A_n.$
\qed

\noindent
Now we turn back to the proof of our proposition.
Let $1\neq \pi\in \Sno$. Then $\pi$ has a representative
$\tilde\pi\in\Sn$ such that for all $i\in I$: $\frac{\tau(i)}{\tau(i)-\fix
\tilde\pi_i}<C$, where $C$ does not depend on $i$. Let $\rho\in\Sno$ be
an arbitrary element. Then it has a representative $\tilde\rho$ such that
all of its component is an even permutation. Thus by Lemma
\ref{lvege2},it is easy to see that $\rho$ is in 
the normal subgroup generated by $\pi$.\,\qed

\begin{theorem}  \label{generic-sofic}
Let $\Gamma$ be a group. The following are equivalent:
\begin{enumerate}
\item	$\Gamma$ is sofic.
\item	There is an injective homomorphism $\Gamma\to\Sno$ for an index set
	$I$ of cardinality at most $|\Gamma|$, and appropriate choice of
	$\tau,\omega$.
\end{enumerate}
In particular, the groups $\Sno$ are themselves sofic for all $I,\tau,\omega$.
\end{theorem}

{\bf Proof of} (1)$\Rightarrow$(2).
For finite $\Gamma$ it is clear.
Let $\Gamma$ be an infinite sofic group.
We define the index set:
$$
I = \left\{ (F,\varepsilon) \ \Big|\
	F\subseteq \Gamma \mbox{ finite, and }
	\varepsilon\in(0,1) \mbox{ rational} \right\}
$$
and for each index $(H,\delta)\in I$ we define the nonempty subset
$$
I_{H,\delta} = \left\{ (F,\varepsilon)\in I \ \Big|\
	H\subseteq F \mbox{ and }\varepsilon\le\delta \right\}
\subseteq I
$$
It is clear that $|I|=|\Gamma|$.
The collection of nonempty subsets $\{ I_{H,\delta}\}$ is closed
under finite intersection, so there is an ultrafilter $\omega$ of
subsets of $I$ containing all $I_{H,\delta}$.
For each index $i=(F,\varepsilon)\in I$ we choose a natural number $\tau(i)$
and a function $\phi_i:\Gamma\to S_{\tau(i)}$ satisfying the
conditions (a)-(c) of Definition~\ref{sofic-group}. Note that by
definition $\lim_{\omega} \epsilon=0$.
By composing $\phi_i$ with the homomorphism
$\coprod_k:S_{\tau(i)}\to S_{k\tau(i)}$ of Definition~\ref{symmetric-group}
(for large enough $k\in\mathbb N$)
we shall achieve that $\lim_\omega \tau(i)=\infty$, and
conditions (a)-(c) of Definition~\ref{sofic-group} remain valid.
The direct product of the functions $\phi_i$ is a function
$\lift\phi:\Gamma\to\Sn$, and composing it with the quotient map
$\Sn\to\Sno$ we get a function $\phi:\Gamma\to\Sno$. We shall prove that
this $\phi$ is an injective homomorphism.
Definition~\ref{sofic-group}(a) implies that
$$
\fix\Big(\phi(e)\phi(f)\,\phi(ef)^{-1}\Big) =
\lim_\omega\fix\Big(\phi_i(e)\phi_i(f)\,\phi_i(ef)^{-1}\Big) \ge
\lim_\omega(1-\varepsilon) = 1
$$
for all $e,f\in \Gamma$.
Therefore $\phi(e)\phi(f)=\phi(ef)$ for all $e,f\in \Gamma$,
so $\phi$ is a group homomorphism.
Definition~\ref{sofic-group}(c) implies that
$$
\fix\phi(e)=
\lim_\omega\fix\phi_i(e) \le
\lim_\omega\varepsilon = 0
$$
for all $1\ne e\in \Gamma$.
But then $\phi(e)\ne 1$ for these $e$,
hence $\phi$ is injective.
\qed

{\bf Proof of} (2)$\Rightarrow$(1).
Let $I,\omega,\tau$ be as in Definition~\ref{sno-def}, we shall prove
that $\Sno$ is sofic. Since subgroups of sofic groups are sofic,
this is enough to prove the Theorem.
So let $F\subset\Sno$ be a finite subset and let $0<\varepsilon<1$ be a number.
For each element $\pi\in\Sno$ we choose a representative
$\lift\pi\in\Sn$, we make sure that $\lift1=1$.
Let $0<\xi<1$ be a real number (to be specified later).
Then for all $\pi,\rho,\sigma\in\Sno$, $\sigma\ne1$ the subsets
$$
\begin{array}{ccl}
A(\pi,\rho) &=& \left\{ i\in I \,\Big|\,
         \fix\Big(\lift\pi(i)\lift\rho(i)\lift{\pi\rho}(i)^{-1}\Big)
         >(1-\xi)\tau(i)\right\}\\
C(\sigma) &=& \left\{ i\in I \,\Big|\,
         \fix\lift\sigma(i) < \frac{\fix\sigma+1}2\tau(i)\right\}
         \kern 40pt (\sigma\ne1)
\end{array}
$$
all belong to $\omega$, so any finite collection of them have nonempty
intersection. Thus we can choose an index $j\in I$ such that
\begin{equation} \label{j-def}
j\in \left(\bigcap_{\pi,\rho\in F}A(\pi,\rho)\right)\bigcap
     \left(\bigcap_{1\ne\sigma\in F}C(\sigma)\right)
\end{equation}
Now we take $n_1=n(j)$, and define the function
$\psi_1:\Sno\to S_{n_1}$ via the formula $\psi_1(\pi)=\lift\pi(j)$.
We define the positive number
$\delta=\max\{\frac{\fix\sigma+1}{2}|1\ne\sigma\in F\}$,
this depends only on $F$ but independent of $\xi$.
Then $\psi_1$ have the following properties:
\begin{enumerate}
\item[(A)]   $\fix\Big(\psi_1(\pi)\psi_1(\rho)\,
             \psi_1(\pi\rho)^{-1}\Big) \ge(1-\xi)n_1$
             \kern20pt for any two elements $\pi,\rho\in F$.
\item[(B)]   $\psi_1(1)=1$.
\item[(C)]   $\fix\psi_1(\sigma)\le\delta n_1$\kern20pt
             for each $1\ne\sigma\in F$.
\end{enumerate}
Now for any natural number $k\in\mathbb N$ we define $n_k=n_1^k$,
and denote by $\psi_k:\Sno\to S_{n_k}$ the composition of $\psi_1$
with the homomorphism $\prod^k:S_{n_1}\to S_{n_k}$ of
Definition~\ref{symmetric-group}.
Then this $\psi_k$ has the following properties:
\begin{enumerate}
\item[(A')]  $\fix\Big(\psi_k(\pi)\psi_k(\rho)\,
             \psi_k(\pi\rho)^{-1}\Big) \ge(1-\xi)^kn_k$
             \kern20pt for any two elements $\pi,\rho\in F$.
\item[(B')]  $\psi_k(1)=1$.
\item[(C')]  $\fix\psi_k(\sigma)\le\delta^k n_k$\kern20pt
             for each $1\ne\sigma\in F$.
\end{enumerate}
We may choose $k$ so large that $\delta^k<\varepsilon$, and then we may
choose sufficiently small $\xi$ to make $(1-\xi)^k>1-\varepsilon$.
With these choices the conditions (A'), (B') and (C') will imply
the conditions (a), (b) and (c) of the Definition~\ref{sofic-group}.
This completes the proof of the theorem.
\qed

%%%%%%%%%%%%%%%%%%%%%%%%%%%%%%%%%%%%%%%%%%%%%%%%%%%%%%%%%%%%%%%%%%%%%%%%

%%  9. oldal

\section{Embedding of countable sofic groups into $R^\omega$}

%%  10. oldal

\begin{theorem} \label{x}
All countable sofic groups are hyperlinear.
\end{theorem}

\proof
According to Theorem~\ref{generic-sofic},
each countable sofic group is isomorphic to a subgroup of one of our
groups $\Sno$ for a countable index set $I$ and appropriate choice
$\tau,\omega$. Hence it is enough to construct a unitary representation
$\phi:\Sno\to R^\omega$ (here we use the same index set $I$ and the
same $\omega$ in the construction of $R^\omega$).

For each index $i\in I$ the symmetric group $S_{\tau(i)}$ has an
irreducible representation
$\sigma_i:S_{\tau(i)}\to\Mat_{\tau(i)\times \tau(i)}$
(unique up to conjugation, in fact the image of a permutation $\pi$
is the linear transformation which permutes the basis vectors via
$\pi$).
Composing it with our embedding $\rho_{\tau(i)}$ we get a unitary
representation $\phi_i:S_{\tau(i)}\to R$. The direct product of these
gives us a unitary representation
$$
\lift\phi = \prod_{i\in I}\phi_i:\Sn\to\prod_{i\in I}R
$$
For each $\pi\in\Sn$ and each index $i\in I$ the matrix
$\sigma_i(\pi(i))$ is a permutation-matrix, so we calculate the
traces:
$\Tr\Big(\sigma_i(\pi(i))\Big) = \fix\pi(i)$ and
$\Tr\Big(1-\sigma_i(\pi(i))\Big)= \tau(i)-\fix\pi(i)$.
Hence we can calculate the normalized traces:
$$
\Tr_R\Big(\phi_i(\pi(i))\Big) = \frac{\fix\pi(i)}{\tau(i)}\le 1
$$
Suppose that $\pi$ and $\rho$ represent
the same element in $\Sno$. Then 
$\lim_\omega \fix(\pi{\rho}^{-1})=1.$ Therefore,
$$\lim_\omega \frac{\Tr(1-\sigma_i(\pi(i)\rho(i)^{-1}))}{\tau(i)}=0,\quad
\lim_\omega \frac{\Tr(1-\sigma_i(\rho(i)\pi(i)^{-1}))}{\tau(i)}=0\,.$$
That is
$$\Tr_\omega ((\lift\phi(\pi)-\lift\phi(\rho))
(\lift\phi(\pi)-\lift\phi(\rho)^*)=0\,.$$
Thus $\lift\phi$ descends to a unitary representation
$\phi:\Sno\to R^\omega$.
Our trace-formula implies, that the image of elements $\pi$ with
$\fix\pi<1$ have $\omega$-trace less than 1, so the representation is
nontrivial. On the other hand $\Sno$ is a simple group according to
Proposition~\ref{sno-properties}.
Hence our $\phi$ must be a faithful representation.
\qed

%%  11. oldal

\begin{corollary}
There exist countable simple non-amenable hyperlinear groups. Also, there
exist finitely generated non-residually amenable hyperlinear groups.
\end{corollary}

\proof
We proved that for a countably index set $I$ and for arbitrary
$\tau,\omega$ as in Definition~\ref{sno-def} our group $\Sno$ is a simple
hyperlinear group. Since each countable sofic group is a subgroup of
one of these, at least one of them must be non-amenable. If $K$ is a 
countable group and $K\subset H$ is a simple group, then there exists
a countable simple group $H'$ , $K\subset H' \subset H$. Hence the corollary
follows. \qed

%%
%%  12. oldal
%%

\section{Characterizations of essentially free amenable actions and
almost-actions}

The goal of this section is to obtain 
Tarski-type characterisations
for essentially free
\\
amenable actions and almost-actions.

\begin{theorem} \label{amenable-action-characterisation}
Let $\Gamma$ be a group acting on a set $X$. For elements
$g\in\Gamma$ we denote by $\Fix(g|X)\subseteq X$ the
fixed point set of $g$ in $X$.
 Then the following
conditions are equivalent:
\begin{enumerate}
\item	The action is an essentially free amenable action.
\item	The action is non-paradoxical in the following sense: $X$
cannot be written as a union of subsets
$X=A_1\cup A_2\cup\dots\cup A_l \cup B_1\cup B_2\cup\dots \cup B_m$
such that
\begin{enumerate}
\item	for some group elements $g_1,h_1,g_2,h_2,\dots,
g_l,h_l\in\Gamma$ the subsets $g_1A_1$, $h_1A_1$, $g_2A_2$, $h_2A_2$,
$\dots$, $g_lA_l$, $h_lA_l$ are pairwise disjoint, 
\item	for each $B_i$ there are elements $p_i\in\Gamma$ with
$B_i=\Fix(p_i| X)$.
\end{enumerate}
%%
%%  13. oldal
%%
\item	For any finite subset $K\subseteq\Gamma$ and any number
$\epsilon>0$ there exists a finite subset $F\subseteq X$ such that
\begin{enumerate}
\item	$\displaystyle\frac{|gF\vartriangle F|}{|F|}<\epsilon$ \ \
	for any $g\in K$, 
\item	$\displaystyle\frac{|\Fix(g | X)\cap F|}{|F|}<\epsilon$ \ \
	for any $g\in K$.
\end{enumerate}
\end{enumerate}
\end{theorem}

\proof
It follows immediately from the next
Theorem~\ref{amenable-almost-action-characterisation}. Or, one can
prove it with the same line of arguments.
\qed

\begin{theorem} \label{amenable-almost-action-characterisation}
Let $\Gamma$ be a group and $X$ be a set. For each $g\in\Gamma$ let
$\phi(g):X\to X$ be a bijection. For elements $p,q\in\Gamma$ we define the
subset
$$
B(p,q)=\{x\in X| \phi(p)\phi(q)\neq \phi(pq)(x)\}
$$
Then the following conditions are
equivalent:
\begin{enumerate}
\item	\label{amenable-cond}
$(\Gamma,X,\mu,\phi)$ is an essentially free amenable
almost-action for some finitely additive probability measure $\mu$.
\item	\label{non-paradoxical-cond}
$\phi$ is non-paradoxical in the following sense:
%%
%%  14. oldal
%%
$X$ cannot be written as a union of subsets
$X=A_1\cup A_2\cup\dots\cup A_l \cup B_1\cup B_2\cup\dots \cup B_m
\cup C_1 \cup C_2\dots \cup C_n $
such that
\begin{enumerate}
\item	for some group elements $g_1,h_1,g_2,h_2,\dots,
g_l,h_l\in\Gamma$ the subsets $g_1A_1$, $h_1A_1$, $g_2A_2$, $h_2A_2$,
$\dots$, $g_lA_l$, $h_lA_l$ are pairwise disjoint,
\item	for each $B_i$ there are elements $p_i, q_i\in\Gamma$ with
$B_i=B(p_i,q_i),$
\item for each $C_i$ there exists an element $r_i\in\Gamma$ with
$C_i=\Fix(\phi(r_i)| X)$\,.
\end{enumerate}
\item	\label{Folner-cond}
For any finite subset $K\subseteq\Gamma$ and any number
$\epsilon>0$
%%
%%  15. oldal
%%
there exists a finite subset $F\subseteq X$ such that
\begin{enumerate}
\item	$\displaystyle\frac{|gF\vartriangle F|}{|F|}<\epsilon$ \ \ \
	for any $g\in K$,
\item	$\displaystyle\frac{|B(p,q)\cap F|}{|F|}<\epsilon$ \ \ \
	for any $p,q\in K$,
\item   $\displaystyle\frac{|\Fix(\phi(g))\cap F|}{|F|}<\epsilon$ \ \ \
        for any $g\in K$.
\end{enumerate}
\end{enumerate}
\end{theorem}

%%
%%  16. oldal
%%

{\bf Proof of
(\ref{amenable-cond})$\Rightarrow$(\ref{non-paradoxical-cond}):}\\
If $(\Gamma,X,\mu,\phi)$ is an essentially free amenable
almost-action, then $\mu(B_i)=0$ and $\mu(C_j)=0$ for all $B_i$ and $C_j$.
Hence $\sum_{i=1}^l\mu(A_i)\geq 1$. Therefore
$$
\sum_{i=1}^l\mu\Big(\phi(g_i)(A_i)\Big) +
\sum_{i=1}^l\mu\Big(\phi(h_i)(A_i)\Big)  \ge 2
$$
providing a contradiction.
\qed

\vskip 7pt
{\bf Proof of
(\ref{Folner-cond})$\Rightarrow$(\ref{amenable-cond}):}\\
We shall denote by $\Fin(X)$ the collection of finite subsets of $X$.
For a fixed pair $K,\epsilon$ denote by
$S_{K,\epsilon}\subseteq\Fin(X)$  the set of finite subsets
$F\subseteq X$ satisfying the conditions of (\ref{Folner-cond}). Then
$S_{K,\epsilon}\cap S_{L,\delta}\supseteq
S_{K\cup L,\min(\epsilon,\delta)}$.
Hence there exists an ultrafilter $\omega$ on the subsets of
$\Fin(X)$ containing each $S_{K,\epsilon}$.
%%
%%  17. oldal
%%
We define a finitely additive measure $\mu$ on $X$ as
$$
\mu(A) = \lim_\omega\frac{|A\cap F|}{|F|}
\kern 25pt\mbox{for all }A\subseteq X.
$$
Clearly $(\Gamma,X,\mu,\phi)$ is an essentially free amenable
almost-action.
\qed

\vskip 7pt
{\bf Proof of
not-(\ref{Folner-cond}) $\Rightarrow$ not-(\ref{non-paradoxical-cond}):}\\
%(\ref{non-paradoxical-cond})$\Rightarrow$(\ref{Folner-cond}):}\\
Let $\epsilon>0$ and $1\in K\subseteq\Gamma$ be a finite subset such
that no $F\subseteq\Gamma$ satisfies the conditions of
(\ref{Folner-cond}). We may suppose that
$1\in K$ and  $K$ is symmetric that is $k\in K$ implies $k^{-1}\in K$.
\begin{lemma}\label{segedlemma}
Let $K^r$ be the set of elements in $\Gamma$ that can be written as the
product of $r$ elements of $K$, then:
$$\bigcup_{g\in K^r} \phi(g)^{-1}((\bigcup_{s,t\in K} B(s,t)\cup
\bigcup_{u\in K} \Fix(\phi(u)| X))\subseteq
\bigcup_{a,b\in K^{r+1}} B(a,b)\cup \bigcup_{c\in K^{2r+1}}
\Fix(\phi(c)|X)\,.$$ \end{lemma}
{\bf Proof}
Suppose that $g\in K^r$, $s,t\in K$ and
$y\in \phi(g)^{-1}(B(s,t))$, that is
$$\phi(s)\phi(t)(\phi(g)(y))\neq \phi(st)\phi(g)(y)\,.$$
Then: $\phi(t)\phi(g)(y)\neq
\phi(tg)(y)$ or
$\phi(s)\phi(tg)(y)\neq\phi(stg)(y)$ or
$\phi(st) \phi(g)(y)\neq \phi(stg)(y)\,.$
Hence
$$y\in\bigcup_{a,b\in K^{r+1}} B(a,b)\,.$$
Now suppose that $g\in K^r, u\in K$ and
$u\in \phi(g)^{-1}(\Fix(\phi(u)|X))\,.$

\noindent
Case 1.: $\phi(g^{-1}ug)(y)\neq\phi(g^{-1}\phi(ug)(y)$ or
$\phi(ug)(y)\neq\phi(u)\phi(g)(y)$ or
$\phi(g^{-1})\phi(g)(y)\neq y$, then:
$$y\in\bigcup_{a,b\in K^{r+1}} B(a,b)\,.$$

\noindent
Case 2.: $\phi(g^{-1}ug)(y)=y\,.$
Then
$$y\in\bigcup_{c\in K^{2r+1}} \Fix(\phi(c)|X)\,.\quad\qed$$

Now let us choose a positive integer $p$ such that $(1+\epsilon)^p>2$.
Define the subsets

%%
%%  18. oldal
%%
%%
%%  19. oldal
%%
$$
A_r = X\setminus(\bigcup_{a,b\in K^{r+1}}B(s,t)\cup 
\bigcup_{c\in K^{2r+1}} \Fix(\phi(c)|X)\,.
$$
By our assumption $A_r$ can not be empty for any $r\geq 1$.
Following the idea of \cite{DSS} consider the bipartite graph $G$ with
vertex sets $(A_p,X)$ such that $(a,x)$ is an edge of $G$ if and only if
$\phi(t)(a)=X$ for some $t\in K^p$.
Observe that for any finite subset $L\subseteq A_p$, the number of neighbours
in $X$ is at least $2|L|$.
Indeed, let $L_0=L$, and define by induction the sequence of subsets
$$
\begin{array}{ccl}
L_s &=& \left\{x\in X\ \Big|\
	\exists g\in K,\ x\in\phi(g)(L_{s-1})\right\}\\
\vspace{-9pt}\\
&=& \left\{x\in X\ \Big|\
	\exists h\in K^s,\ x\in\phi(h)(L)\right\}\\
\end{array}
\kern35pt 1\le s\le p
$$

It is clear from our lemma that $L_s\cap
(\bigcup_{p,q\in K} B(p,q)\cup
\bigcup_{r\in K} \Fix(\phi(r)| X)=\emptyset\,$.
 Hence
$|L_s|\ge(1+\epsilon)|L_{s-1}|\ge(1+\epsilon)^s|L|$ for all $s$.
The set of neighbours (in $G$) of the elements of $L$ is just $L_p$
with at least $2|L|$ elements.
Therefore we can apply Hall's (2,1)-matching theorem to our bipartite
graph $G$:
there exist matchings $m_1:A\to X$ and $m_2:A\to X$ whose images are
disjoint. If $s\ne t\in K^p$, let
$$
Q_{s,t}=\left\{a\in A_p\ \Big|\
	m_1(a)=\phi(s)(a),\ m_2(a)=\phi(t)(a) \right\}
\ .
$$
Then
$$
A_p = \bigcup_{s\ne t\in K^p} Q_{s,t}
$$
and
$$
\Bigg\{\phi(s)\Big(Q_{s,t}\Big),\
	\phi(t)\Big(Q_{s,t}\Big)\Bigg\}_{s\ne t\in K^p}
$$
are disjoint subsets. Hence $X$ can be covered the way described in
(\ref{non-paradoxical-cond}).
\qed

%%
%%  20. oldal
%%

The corollary below easily follows from (\ref{Folner-cond}) of
Theorem~\ref{amenable-almost-action-characterisation}.

\begin{corollary}
Let $\Gamma$ be a discrete group. Then $\Gamma$ has an essentially
free almost-action if and only if $\Gamma$ is sofic.
\end{corollary}

For essentially free amenable actions we have the following
proposition.

\begin{proposition}
If $\Gamma$ is amenable or residually finite then it has an
essentially free amen\-able action.
\end{proposition}

\proof
If $\Gamma$ is amenable then the left action on itself is a free
amenable action. If $\Gamma$ is residually finite, then let $X$ be the
disjoint union of all finite quotients of $\Gamma$ with the natural
left $\Gamma$-actions. The finitely additive measure is defined with
an ultra-limit, as in the proof of
Theorem~\ref{amenable-almost-action-characterisation}.
\qed

%%
%%  21. oldal
%%

Now, we show the existence of groups which have no essentially free
amenable actions.

\begin{proposition}
Let $\Gamma$ be a discrete group of Kazhdan's Property (T). Then
$\Gamma$ has an essentially free amenable action if and only if it is
residually finite.
\end{proposition}

\proof
Suppose that $(\Gamma, X, \mu)$ is an essentially free amenable action
of a Kazhdan group $\Gamma$. Let $X_F\subset X$ be the union of the
finite $\Gamma$-orbits. First suppose that $\mu(X_F)=1$. Then for any
$1\ne g\in\Gamma$ there exists a finite orbit on which $g$ acts
non-trivially. Otherwise, the fixed point set of $g$ had positive measure.
Hence $\Gamma$ is residually finite in this case.
Now suppose that $\mu(X_F)<1$.
%%
%%  22. oldal
%%
Then $\mu(Y)>0$ where $Y$ is the union of the infinite orbits. Thus
$\hat\mu=\frac{\mu}{\mu(Y)}$ defines a finitely additive probability
measure on $Y$.
Consider the Hilbert space $\mathcal H$ of the square-summable
functions on $Y$. Then $\Gamma$ acts on $\mathcal H$ without a fixed
unit vector. On the other hand, by
Theorem~\ref{amenable-action-characterisation}~(\ref{Folner-cond})
for any $\delta>0$ and all finite $K\subseteq \Gamma$
there exists a finite set $F\subseteq Y$ such that
$$
\left\|	\frac{\chi_F}{\sqrt{|F|}}-g\left(\frac{\chi_F}{\sqrt{|F|}}\right)
	\right\|_{\mathcal H} < \delta
\kern 25pt\mbox{for all }g\in K\ ,
$$
where  $\chi_F$ denotes the characteristic function of $F$. Since
$\Gamma$ has Kazhdan's Property (T), it must have a fixed unit vector
in $\mathcal H$, leading to a contradiction.
\qed

\noindent
We do not know any example of a non-residually finite sofic group of
Kazhdan's Property (T). However we have the following proposition.
\begin{proposition} \label{kazhdan}
If all sofic groups of Kazhdan's Property (T) is residually finite, then
there exists a non-residually finite word-hyperbolic group.
\end{proposition}
\begin{lemma}
If $\Gamma\in\Sigma_m$ (with $m$ specified generators) is a limit
point of sofic elements, then $\Gamma$ is sofic as well.
\end{lemma}

\proof
Let $\lim_{n\to\infty}\Gamma_n=\Gamma$ (with $m$ specified generators)
in $\Sigma_m$.
Suppose that $F\subseteq\Gamma$ is a finite set. Then there exists an
index $k$ and an embedding $t_F:F\to\Gamma_k$ such that
$t_F(ab)=t_F(a)t_F(b)$ whenever $a,b,ab\in F$.
Since soficity is a property decided by the properties of
multiplication restricted to such finite subsets $F$, the soficity of
all $\Gamma_k$ implies the soficity of $\Gamma$.
\qed

%%
%%  31. oldal
%%
\noindent
{\bf Proof of the proposition:}
Champetier \cite{Champ} showed that the closure of hyperbolic groups in
$\Sigma_m$ contains Kazhdan groups that have no subgroups of finite
index. Thus if all hyperbolic group were residually finite, then there
exists a sofic group of Kazhdan's Property (T).
\qed

%%
%%  23. oldal
%%

\section{Obstruction for essential freeness and inner amen\-ability}

Let $(\Gamma,X,\mu)$ be an action of a group $\Gamma$ on a
set $X$ with a $\Gamma$-invariant finitely additive measure $\mu$. We
denote the class of such actions by $\Inv_\Gamma$.
Let $N_{\Gamma,X,\mu}$ be the set of elements $g\in\Gamma$ such that
$\mu(\Fix(g| X))=1$.
Then $N_{\Gamma,X,\mu}$ is a normal subgroup of $\Gamma$.
Indeed,
$\Fix(g| X)\cap\Fix(h| X)
\subseteq\Fix(gh| X)$
implies that
$\mu(\Fix(gh| X))=1$ if
$\mu(\Fix(g| X))=\mu(\Fix(h| X))=1$.
Also, $\Fix(aga^{-1}| X)=a\Fix(g| X)$.
Let $N_\Gamma = \bigcap_{(\Gamma,X\mu)\in\Inv_\Gamma} N_{\Gamma,X,\mu}$,
this is a characteristic subgroup in $\Gamma$.
%%
%%  24. oldal
%%  25. oldal
%%  26. oldal
%%  27. oldal
%%
\begin{proposition} \label{N-Gamma}
Let $\Gamma$ be a group.
\begin{enumerate}
\item	\label{0-1-action-for-N-Gamma}
There is an action $(\Gamma,X,\mu)\in\Inv_\Gamma$ such that
$
\mu(\Fix(g| X))=\left\{
	\begin{array}{cl}
	1&\mbox{for }g\in N_\Gamma\\
	0&\mbox{for }g\in\Gamma\setminus N_\Gamma\\
	\end{array}
\right.
$
\item	\label{N-Gamma=1}
$\Gamma$ has an essentially free amenable action on a set $X$
if and only if $N_\Gamma =1$.
\item	\label{Gamma/N-Gamma-sofic}
$\Gamma/N_\Gamma$ is a sofic group.
\end{enumerate}
\end{proposition}

{\bf Proof of (\ref{0-1-action-for-N-Gamma}):}
For each element $g\in\Gamma\setminus N_\Gamma$ we choose an action
$(\Gamma,X_g,\mu_g)\in\Inv_\Gamma$ such that
$\mu_g(\Fix(g| X_g))<1$.

We define the index set $I$ of all pairs $(F,\varepsilon)$ where
$F\subseteq\Gamma\setminus N_\Gamma$ is a finite subset and
$0<\varepsilon<1$ is a real number.
As in the proof of Theorem~\ref{generic-sofic}, there is an
ultrafilter $\omega$ on the subsets of $I$ which contains all subsets
of the form
$
I_{H,\delta} =
\{ (F,\varepsilon)\in I\ |\ H\subset F,\ \varepsilon<\delta\}
$
with some $(H,\delta)\in I$.
For each index $(F,\varepsilon)$ there is an integer $K_{F,\varepsilon}$ such
that
$$
\mu_g(\Fix(g| X_g))^{K_{(F,\varepsilon})}<\varepsilon
\kern 25pt \mbox{for all }g\in F\,.
$$
If $(X,\mu_X),(Y,\mu_Y)\in \Inv_\Gamma$ then $(X\times Y,\mu_X\times\,\mu_Y)$
is defined the following way:
$$\mu_X\times\mu_Y(H)=\int H_x \,d\mu_X\,,$$
where $H_x=\{y\in Y| (x,y)\in H\}$ and
$\int\,d\mu_X$ is the translation invariant functional on $L^\infty(X)$ 
associated to $\mu_X$ \cite{Pat}.
Inductively, one can define
$(\prod^n_{i=1} X_i,\prod^n_{i=1}\mu_{X_i})$ as well. For a fixed
$(F,\epsilon)\in I$, let
$$(X_{(F,\epsilon)},\mu_{(F,\epsilon)})=\prod_{g\in F}(\prod^{K_{(F,\epsilon)}}
_{i=1} X^i_g)\,,$$
where $X^i_g$ is a copy of $X_g$. Then for any $g\in F$,
$\mu(\Fix(g| X_{(F,\epsilon)})\leq \epsilon\,.$
Let $X=\coprod_{(F,\epsilon)} X_{(F,\epsilon)}$ and
$$\mu(H)=\lim_\omega\,\mu_{(F,\epsilon)}(H\cap X_{(F,\epsilon)})\,.$$
Then $(X,\mu)\in\Inv_\Gamma$ and $\mu(\Fix(g|X))=\lim_\omega \epsilon=0\,,$
if $g\in \Gamma\backslash N_\Gamma$.
This proves our claim.
\qed

{\bf Proof of (\ref{N-Gamma=1}):}
If $(\Gamma,X,\mu)$ is an essentially free amenable action, then
$(\Gamma,X,\mu)\in\Inv_\Gamma$, and $N_{\Gamma,X,\mu}=1$.
Hence $N_\Gamma=1$.
On the other hand, if $N_\Gamma=1$ for a group $\Gamma$, then in
(\ref{0-1-action-for-N-Gamma}) we build an action
which is clearly essentially free and amenable.
\qed

{\bf Proof of (\ref{Gamma/N-Gamma-sofic}):}
Let $(\Gamma,X,\mu)\in\Inv_\Gamma$ be the action constricted in
(\ref{0-1-action-for-N-Gamma}).
Then the group $\Gamma/N_{\Gamma,X,\mu}=\Gamma/N_\Gamma$
has an essentially free amenable almost-action on $(X,\mu)$:
each element $g$ acts via an arbitrarily choosen lifting
$\tilde g\in\Gamma$. Two different liftings act identically
outside a subset of $\mu$-measure 0, hence we obtain an
almost action. Essential freeness follows from the fact that
$\mu(\Fix(\tilde g|X))=0$ unless $g=1$.
\qed

%%
%%  28. oldal
%%
\begin{corollary} \label{inner}
An inner-amenable discrete group either has a nontrivial sofic
quotient, or it is almost commutative (\cite{BH}).
\end{corollary}

\proof
Let $\Gamma$ be an inner amenable group, i.e. its action on
$X=\Gamma\setminus\{1\}$ is amenable.
Let $\mu$ be a conjugation invariant finitely additive
probability measure on $X$. If $N_{\Gamma,X,\mu}\ne\Gamma$,
then by Proposition~\ref{N-Gamma} the nontrivial quotient
$\Gamma/N_\Gamma$ is sofic.
On the other hand, if $N_{\Gamma,x,\mu}=\Gamma$, then for any finite
set $\{a_1,a_2,\dots a_k\}\subseteq\Gamma$:
$$
\mu\Big(\Fix(a_1|X)\cap\Fix(a_2|X)\cap\dots\cap\Fix(a_k|X)\Big)=1
\ .
$$
Hence there is an element $b\in X=\Gamma\setminus\{1\}$ which is fixed
by each $a_i$. Then $b$ commutes with $a_1,a_2,\dots a_k$,
that is $\Gamma$ is almost commutative.
\qed

%%
%%  29. oldal
%%

\section{Approximating the $L^2$-invariants}
Let $\Gamma$ be a finitely generated sofic group and $(V_n,E_n)_{n\geq 1}$ be a sequence of
finite graphs as in the Definition \ref{sofic1}.
 Let $A=\{a_{ij}\}_{1\leq i,j \leq d}\in\Mat_{d\times d}
(\bZ \Gamma)$ be a positive, self-adjoint operator. Consider the operator kernel 
of $A$, that is the function $K_A:\Gamma\times \Gamma\to \Mat_{d\times d}(\bZ)$ such
that if $f\in [l^2(\Gamma)]^d$ a vector-valued function then:
$$Af(x)=\sum_{y\in \Gamma} K_A(x,y) f(y)\,.$$
Obviously, $K_A(x,y)=A^\gamma$ if $A=\sum_{\gamma\in \Gamma}
A^\gamma\cdot \gamma$, $A^\gamma=\{a_{ij}^\gamma\}_{1\leq i,j \leq d}
\in \Mat_{d\times d}(\bZ)$
and $x=\gamma y$.
Note that there exists a constant $w_A$, the width of $A$, such that 
$K_A(x,y)=0$ if $d(x,y)>w_A$ in the word-metric of $\Gamma$ with respect to
the generating system $S$.

\noindent
Now we construct the approximation kernel, for $m>n_{(w_A,\frac{1}{2})}$,
$K^m_A:V_m\times V_m\to \Mat_{d\times d}(\bZ)$.
Let 
$K^m_A(x,y)=0$, if $y\in V^0_m$ and
let $K^m_A(x,y)=K_A(\gamma,1_\Gamma)$ if $y\in V^0_m, x=\psi_y(\gamma)$.
(We use the notation of Definition \ref{sofic1}.)

\noindent
We denote by $A_m$ the bounded linear transformation on $[l^2(V_m)]^d$
defined by the kernel function $K^m_A$. The goal of this section is to 
prove the following proposition:
\begin{proposition} \label{schick}
\begin{description}
\item{(a)} $\lim_{m\to\infty} \frac{\dim \ker (\Delta_m)}{|V_m|}=F_A(0)\,.$
\item{(b)} $\limsup_{m\to\infty} \frac{\ln \det^* (\Delta_m)}{|V_m|}\leq
\ln \det_{\cN}(\Delta)\,,$  
\end{description}
where $\Delta=A^*A$, $\Delta_m=A^*_m A_m$ and $\det^*(\Delta_m)$ is
the product of the non-zero eigenvalues of $\Delta_m$.
\end{proposition}
The proposition is the sofic version of the approximation theorem of
L\"uck, first proved for residually finite groups \cite{Lueck} and
later extended for residually amenable groups by Clair \cite{Clair}
and Schick \cite{Schick1}.
\begin{lemma}\label{lemma l10}
There exists a constant $K$ such that $K\geq \|\Delta\| $ and $K\geq
\|\Delta_m\|$ for all $m\geq 1$, where $\|\,\|$ denotes the operator norm.
\end{lemma}
\proof Let $B=\{b_{ij}\}_{1\leq i,j\leq n}
\subset\Mat_{n\times n} (\bR)$ be a matrix
with at most $L$ non-zero entries in each of its rows and columns.
Let $M$ be $\sup_{1\leq i,j\leq n} |b_{ij}|$. Then if we regard
$B$ as an operator on the Euclidean space $\bR^n$,
\begin{equation}
\label{10e}
\|B\|\leq LM \,.
\end{equation}
Indeed, let $v=\sum^n_{i=1} c_i d^i\in \bR^n$, then $\|v\|^2=\sum^n_{i=1}
c_i^2\,.$
$$\|Bv\|^2=\sum^n_{i=1}(\sum^n_{j=1} b_{ij}c_j)^2\leq M^2
(\sum^n_{i=1}\chi_{ij}c_j)^2\,$$
where $\chi_{ij}=1$ if $b_{ij}\neq 0$ and $\chi_{ij}=0$, if $b_{ij}=0\,.$
Since any row contains at most $L$ elements,
$$(\sum^n_{j=1}\chi_{ij}c_j)^2\leq L(\sum^n_{j=1}\chi_{ij}c^2_j)\,. $$
Also, each column contains at most $L$ elements, hence
$$\sum^n_{i=1} L(\sum^n_{j=1}\chi_{ij}c_j^2)\leq L^2(\sum^n_{j=1} c^2_j)\,.$$
Thus, $\|Bv\|^2\leq L^2M^2\|v\|^2$ and (\ref{10e}) follows. It is easy to see
that (\ref{10e}) implies the statement of our Lemma.\qed
\begin{lemma}\label{l11}
For any polynomial $p(x)$,
$$\Tr_\Gamma(P(\Delta))=\lim_{m\to\infty}\frac{\Tr(p(\Delta_m))}{|V_m|}\,.$$
\end{lemma}
\proof
Let $p$ be a polynomial of rank $s$. Then for large $m$,
$$K^{P(\Delta_m)}(x,y)=K^{P(\Delta)}(\gamma,1_\Gamma)\,,$$
if $|\gamma|\leq 2s w_A, y\in V^0_m$ and $K^{P(\Delta)}$ resp.
$K^{P(\Delta_m)}$ are the operator kernels of $P(\Delta)$ resp. $P(\Delta_m)$.
$$\frac{\Tr P(\Delta_m)}{|V_m|}=\frac{\sum_{z\in V_m} \tra K^{P(\Delta_m)} (z,z)}{|V_m|}, $$
$$Tr_\Gamma P(\Delta)= \tra \,K^{P(\Delta)}(1_\Gamma,1_\Gamma)\,.$$
By the uniform boundedness of the vertex degrees of $(V_m,E_m)$ there is
a uniform bound $K$ on the coefficients of the matrices
$K^{P(\Delta_m)}(x,y)$,
depending only on $p$ and $A$.
Hence:
$$|V^0_m|\tra \,K^{P(\Delta)}(1_\Gamma,1_\Gamma)-d\cdot K \cdot  |V_m\setminus
V^0_m|\leq $$
$$
\sum_{z\in V_m}\tra K^{P(\Delta_m)}(z,z)\leq |V^0_m|
 \tra\, K^{P(\Delta)}(1_\Gamma,1_\Gamma)+d\cdot K\cdot
|V_m\setminus V^0_m|\,. $$
Since $\frac{|V^0_m|}{|V_m|}\to 1$, the Lemma follows. \qed
\begin{lemma}\label{l13}
$\ln \det^*(\Delta_m)\geq 0$, for any $m>0$.\end{lemma}
\proof
The product of the non-zero eigenvalues of $\Delta_m, \,\det^*(\Delta_m)$
is the lowest non-zero coefficient of the characteristic polynomial
of $\Delta_m$. Since $\Delta_m$ has only integer coefficients, this number is
a positive integer.
\qed

\noindent
Now Proposition \ref{schick} follows exactly the same way as Theorem 6.9
\cite{Schick1}. \qed
\begin{theorem}\label{soficdet}
The Determinant Conjecture holds if $\Gamma$ is a sofic group.
\end{theorem}
\begin{theorem}\label{soficlimit}
If $(\Gamma_n,S_n)_{n\geq 1}$ is a convergent sequence of torsion-free sofic groups for which
the Atiyah Conjecture holds, then it holds for the limit group as
well.
\end{theorem}
First of all, Proposition \ref{schick} and Lemma \ref{l13} immediately show
that if $\Delta$ is a positive self-adjoint operator, then
$\ln \det (\Delta^2)\geq 0$. By \cite{Lueck2} $\ln \det (\Delta^2)=
2 \ln \det (\Delta)$. Hence, Theorem \ref{soficdet} follows.
Now, let $(\Gamma_n,S_n)_{n\geq 1}$ be marked sofic groups of $m$ generators such
that $\lim_{n\to\infty} (\Gamma_n,S_n) = (\Gamma,S)$ in $\Sigma_m$.
We may suppose that $B_{\Gamma_n}(l)\simeq B_\Gamma(l)$ for all $n$ and that $B_\Gamma(l)$
contains the support of $\Delta\in \Mat_{d\times d} (\bZ \Gamma)\,.$
For each $\Gamma_n$, let us choose an approximating sequence of finite directed 
graphs $\{V_{p,n},E_{p,n}\}_{p\geq 1}$ as in the definition of sofic groups.
The following lemma is easy to prove and is left for the reader.
\begin{lemma} \label{l14}
For any integer sequence $\{p_i\}_{i\geq 1}$ one can pick
directed graphs $\{V_{r_i,i}, E_{r_i,i}\}_{i\geq 1}$ from the graphs above
so that they are approximating $\Gamma$ and $r_i\geq p_i$ for all $i\geq 1$.
\end{lemma}
\qed

\noindent
Since $B_{\Gamma_n}(l)\simeq B_\Gamma(l)$, one can pull back the operator kernel of $\Delta$ to $\Gamma_n$
obtaining the operators $\Delta_n\in \Mat_{d\times d} (\bZ \Gamma_n)$.
Also, for each $n$, one can consider the approximating operators
$\Delta_{p,n}$ associated to the finite graphs $\{V_{p,n},E_{p,n}\}$.
Then by Proposition \ref{schick}
$$\dim_{\Gamma_n} \ker (\Delta^*_n\Delta_n)=\lim_{p\to \infty} 
\frac{\ker (\Delta^*_{p,n}\Delta_{p,n})}{|V_{p,n}|}\,.$$
Now let us suppose that
$$ 
|\dim_{\Gamma_n} \ker (\Delta^*_n\Delta_n)-
\frac{\ker (\Delta^*_{p,n}\Delta_{p,n})}{|V_{p,n}|}|<\frac{1}{n}\,$$
if $p\geq p_n\,.$
Pick a sequence $\{V_{r_n,n}, E_{r_n,n}\}_{n\geq 1}$, $ r_n\geq p_n$
as in the previous lemma.
Then
$$\dim_{\Gamma} \ker (\Delta^*\Delta)=\lim_{n\to \infty}
\frac{\ker (\Delta^*_{r_n,n}\Delta_{r_n,n})}{|V_{r_n,n}|}\,.$$
Thus
$$\dim_{\Gamma} \ker (\Delta^*\Delta)=\lim_{n\to \infty}
\dim_{\Gamma_n} \ker (\Delta^*_n\Delta_n)\,. $$
Since $\ker 
(\Delta^*_n,\Delta_n)=\ker(\Delta_n)$, Theorem \ref{soficlimit}
follows.\qed
\begin{corollary}
The Determinant Conjecture holds for all LEF-groups, for some
countable non-amenable simple groups and for some finitely generated
non-residually amenable groups as well.
\end{corollary}


\begin{thebibliography}{9}
\bibitem{BH} {\sc E. Bedos and P. de la Harpe},
Moyennabilit\'e int\'erieure des groupes: d\'efinitions et exemples,
{\sl Enseign. Math. (2)} {\bf 32} (1986) no. 1-2, 139 - 157
\bibitem{Brenner} {\sc J. L. Brenner},
Covering theorems for FINASIGs. VIII.
Almost all conjugacy classes in $A_n$ have exponent $\leq 4$.
{\sl J. Austral. Math. Soc. Ser. A.} {\bf 25} (1978) no. 2, 210-214.
\bibitem{Champ} {\sc C. Champetier},
L'espace des groupes de type fini,
{\sl Topology} {\bf 39} (2000) no. 4, 657-680.
\bibitem{Clair} {\sc B. Clair},
Residual amenability and the approximation of $L^2$-invariants,
{\sl Michigan Math. J.} {\bf 46} (1999) no. 2, 331-346.
\bibitem{DSS}{\sc W. A. Deuber, M. Simonovits and V. T. S\'os},
A note on paradoxical metric spaces,
{\sl Studia Sci.Hung.Math.} {\bf 30} (1995), no. 1-2, 17-23.
\bibitem{ESZ1} {\sc G. Elek and E. Szab\'o},
On sofic groups, {\sl to appear in the Journal of Group Theory}
(http://arXiv.org/abs/math/0305352).
\bibitem{ESZ2} {\sc G. Elek and E. Szab\'o}
Sofic groups and direct finiteness, {\sl to appear in the Journal of Algebra}
(http://arXiv.org/abs/math.RA/0305440).
\bibitem{Gro} {\sc M. Gromov},
Endomorphisms of symbolic algebraic varieties,
{\sl J. Eur. Math. Soc.} {\bf 1} (1999) no. 2, 109-197.
\bibitem{Kirch} {\sc E. Kirchberg}, On nonsemisplit extensions, tensor
products and exactness of group $C^\star$-algebras,
{\sl Invent. Math.} {\bf 112} (1993) no. 3, 449-489 .
\bibitem{Lueck} {\sc W. L\"uck},
Approximating $L^2$-invariants by their finite-dimensional analogues,
{\sl Geom. Funct. Analysis} {\bf 4} (1994) no. 4, 455-481.
\bibitem{Lueck2} {\sc W. L\"uck},
$L^2$-invariants: theory and applications to geometry and $K$-theory. 
{\sl Ergebnisse der Mathematik und ihrer Grenzgebiete. 3.}
{\bf 44} Springer Verlag (2002)
\bibitem{Oza} {\sc N. Ozawa}, About the QWEP conjecture,
{\sl International Journal of Math.} {\bf 15} (2004) no. 5, 501-530.
\bibitem{Pat} {\sc A. L. T. Paterson},
Amenability,
{\sl Mathematical Surveys and Monographs} {\bf 29} {\it
American Mathematical Society, Providence} (1988)
\bibitem{Rad} {\sc F. Radulescu},
 The von Neumann algebra of the non-residually
 finite Baumslag group \\ $< a,b | a b^3 a^{-1} = b^2 >$ \\ embeds
into $ R^\omega$,
{\sl preprint 2002}
(http://www.arxiv.org/abs/math.OA/0004172).
\bibitem{Schick1} {\sc T. Schick},
$L^2$-determinant class and approximation of $L^2$-Betti numbers
{\sl Trans. Amer. Math. Soc.} {\bf 353} (2001) no. 8, 3247-3265.
\bibitem{Schick2} {\sc T. Schick},
Integrality of $L^2$-Betti numbers,
{\sl Math. Ann.} {\bf 317} (2000) no. 4, 727-750.
\bibitem{VG} {\sc A. M. Vershik and E. I. Gordon},
Groups that are locally embeddable in the class of finite groups,
{\sl Algebra i Analiz} {\bf 9} (1997) no. 1 71-97.
\bibitem{Wei} {\sc B. Weiss},
Sofic groups and dynamical systems, ({\it Ergodic theory and harmonic analysis,
Mumbai, 1999)}
{\sl Sankhya Ser. A.} {\bf 62} (2000) no. 3,  350-359.


\end{thebibliography}
\end{document}